\documentstyle[12pt]{article}
\newcommand{\const}{\mathop{\rm const}\limits}

\newcommand{\vraisup}{\mathop{\rm vraisup}\limits}

\newcommand{\mod}{\mathop{\rm mod}\limits}

\begin{document}

\begin{center}

{\bf TAIL ESTIMATES FOR MARTINGALE UNDER "LLN" NORMING SEQUENCE }\\

\vspace{4mm}

{\bf E. Ostrovsky, L.Sirota.}\\

\vspace{3mm}

Department of Mathematics, Bar-Ilan University, Ramat-Gan, 59200,Israel.\\
e-mail: \ galo@list.ru, eugostrovsky@list.ru \\

Department of Mathematics, Bar-Ilan University, Ramat-Gan, 59200,Israel.\\
e - mail: \ sirota@zahav.net.il\\

\end{center}

{\bf Abstract.} \par

\vspace{3mm}

  In this paper non-asymptotic exponential and moment
estimates are derived for  tail of distribution for discrete time martingale
under norming sequence 1/n, as in the classical Law of Large Numbers (LLN),
by means of martingale differences as a rule in the terms of {\it unconditional}
moments and tails of distributions of summands. \par
 We show also the exactness of obtained estimations. \par

\vspace{3mm}

{\it Key words:} Random variables and vectors,
martingales, martingale differences, Law of Large Numbers (LLN),
lower and upper estimates, great or large deviations, moment, Banach spaces of
random variables, tails 
of distribution, conditional expectation. \par

\vspace{3mm}
 {\it Mathematics Subject Classification (2002):} primary 60G17; \ secondary
60E07; 60G70.\\

\vspace{4mm}
\begin{center}
\section{ Introduction. Notations. Statement of problem.} \par
\end{center}
\vspace{3mm}

 Let $ (\Omega,F,{\bf P} ) $ be a probabilistic space,
$ \xi(1), \xi(2), \ldots,\xi(n), \ n < \infty $ being a  centered
$ ({\bf E} \xi(i) = 0, i=1,2,\ldots,n) $ martingale - differences on the basis of the
{\it flow }of $ \sigma - $ fields (filtration) $ F(i): F(0) =
\{\emptyset, \Omega \}, \ F(i) \subset F(i+1) \subset F, \
\xi(0) = 0; \ {\bf E} |\xi(i)| < \infty, $ and for every $ i \ge 0,
 \ \forall k = 0,1,\ldots,i -1 \ \Rightarrow  $

$$
{\bf E} \xi(i)/F(k) = 0; \ {\bf E}\xi(i)/F(i) = \xi(i) \  (\mod \ {\bf P}).
$$

\vspace{3mm}
 We denote

 $$
 S(n) = \sum_{i=1}^n \xi(i),  \ {\bf Q}_n(x) = {\bf P} (S(n)/n > x), \ x > 0.
 $$
  Obviously, if the sequence $ \{\xi(i) \}  $  satisfies the Weak Law of Large Numbers (WLLN),
 for instance, the r.v. $ \{\xi(i) \} $ are i.i,d. with finite expectation $ {\bf E} |\xi(1)| < \infty, $
 then

 $$
 \lim_{n \to \infty} {\bf Q}_n(x)=0.
 $$

\vspace{4mm}

{\bf Our aim is to obtain the tail estimates $ {\bf Q}_n(x) $ for  $ S(n)/n $  via the
moment and tail estimates (conditional or not) of the sequences } $  \{ \xi(i) \}. $ \par

\vspace{3mm}

 Our estimates improve or generalize the well-known inequalities belonging to   Fan X., Grama I, and Liu Q
\cite{Fan1}, Grama I.G. \cite{Grama1}, Grama I. and Haeusler E. \cite{Grama2},  Lesign E. and Volny D.
\cite{Lesign1}, Li Y.  \cite{Li1}, Liu Q. and Watbled F. \cite{Liu1} etc.\par

 Notice that in the articles \cite{Lesign1}, \cite{Liu1}, \cite{Volny1} and many others are described some
{\it new} applications of these estimates in the theory of dynamical system  and in the theory
of polymers. \par
 Sometimes, see e.g.  \cite{Lesign1}, \cite{Li1}, \cite{Liu1} etc. the non-asymptotical
estimates of probability $ {\bf Q}_n(x /\sqrt{n})  $ are called "large deviation" (or "great deviation"). \par

  By our opinion, there exists a duality in this terminology:  in many books and articles
\cite{Ibragimov1},  \cite{Petrov1}, \cite{Saulis1},  \cite{Grama1}, \cite{Grama2}, \cite{Ratchkauskas1}
and many others the notion "large deviation" (or "great deviation") was used for the asymptotical
behavior of the fraction

$$
z := \frac{{\bf Q}_n(x/\sqrt{n})}{1-\Phi(x)}, \ \Phi(x) = (2 \pi)^{-1/2} \int_{-\infty}^x \exp(-y^2/2) \ dy
$$
as $ x = x(n) \to \infty. $ \par

 We denote as usually the $ L(p) $ norm of the r.v. $ \eta $ as follows:

 $$
 |\eta|_p = \left[ {\bf E} |\eta|^p  \right]^{1/p}.
 $$

\vspace{4mm}

 The paper is organized as follows.  The second section contains the  main result: estimation of tail-function for normed sum
of martingale differences.  In the third section  we intend to show the exactness of our estimates.\par
 Fourth section contains is devoted to the tail estimates for martingales under more strict conditions;
we obtain at the same estimates as in the independent case; following, they are not improvable. \par
In the last section we consider an inverse statement  to our problem.\par
 In all the sections we bring some examples in order to show the exactness of obtained estimates. \par

\vspace{4mm}

\section{ Main result: estimation of tail-function for normed sum. }

\vspace{4mm}
 We need for beginning to introduce some notations. Let  $ T(x) $ and $ G(x), \ x > 0 $ be a two
tail - functions, i.e.  such that $ T(0) = G(0)=1, \ T(\cdot), G(\cdot) $ are monotonically
decreasing, right continuous and $ T(\infty) = G(\infty) = 0. $ Let us denote
$$
T \vee G(x) =  \min(4 \inf_{y > 0} (T(y)  + G(x/y) ),1).
$$
 The function $ T \vee G(x) $ has a following sense: if $ T(\xi,x) \le T(x), \
T(\eta,x) \le G(x), $ then
$$
 T(\xi \cdot \eta,x) \le T \vee G(x).
$$

  For example, if $ \forall x > 0 \ T(\xi,x)\le \exp\left(-x^{q(1)} \right), \ T(\eta,x) \le \exp
 \left(-x^{q(2)} \right), \  q(1),q(2) = \const > 0, $  then
$$
T(\xi \cdot \eta,x) \le \min \left( 1, 4 \exp(-C(q(1),q(2)) \ x^{q(1)q(2)/(q(1)+q(2)}
\right).
$$
 Let $ \xi $  be a random variable; its {\it tail-function } $ T(\xi,x), x \ge 0 $  we introduce
as follows:

$$
T(\xi,x) = \max({\bf P}(\xi \ge x), {\bf P}(\xi \le - x)).
$$

 Further, for the tail - function $ T(\cdot) $  we denote the following operator
(non - linear)
$$
W[T](x) = \min \left(1, \inf_{v > 0} \left[ \exp(-x^2/(8v^2)) - \int_v^{\infty}
x^2 \ dT(x) \right] \right),
$$
if there exists the second moment
$$
\int_0^{\infty} x^2 \ |dT(x)| < \infty.
$$

\vspace{3mm}

{\bf Theorem 2.1.} Let again $ \xi(i)  $  be a sequence of martingale-differences relative to
some filtration $ \{F(i)\} $ and $ T(\xi(i),x) \le T(x), \ T(x) $ be some tail-function. Then at $ x \ge 2 $

$$
Q_n(x) \le  W[T]( x \sqrt{n}). \eqno(2.1)
$$

\vspace{3mm}
{\bf Proof.}  We will use the following fact, see \cite{Ostrovsky4}, Lemma 1. \par
Let again $ \xi(i)  $  be a sequence of martingale-differences relative to
some filtration $ \{F(i)\} $ and $ T(\xi(i),x) \le T(x), \ T(x) $ be some tail-function. Then at $ x \ge 2 $

$$
\sup_{b: \sum b^2(i) = 1} T \left(\sum_i b(i) \xi(i), x \right)  \le W[T](x). \eqno(2.2)
$$

 We obtain the proposition (2.1) after choosing $ b(i) = 1/\sqrt{n} $ and substituting $ x := x \sqrt{n}. $\par

\vspace{3mm}

{\bf Example 2.1.} Assume that $ T(\xi(i),x) \le Y \exp(-(x/K)^q), \ K,x,q > 0,
 Y = \const \ge 1. $ Denote here
$$
\delta = \delta(q) = (\min(q/2,1))^{-1/q};
$$
$$
q \in (0,2] \Rightarrow \beta(q) = \max \left(\frac{1}{q} \Gamma \left( \frac{2}{q}\right),
\frac{e}{q} \left(\frac{2}{eq} \right)  \right);
$$
$$
q > 2 \Rightarrow \beta(q) = \sup_{v \ge 0} \exp \left(v^q \right)\int_v^{\infty} x
\exp \left(-x^q \right) \ dx  \le \Gamma(2/q) /(qe);
$$
 $ \Gamma(\cdot) $ is usually Gamma - function. We obtain after some
calculations that at $ x \ge 2 $
$$
Q_n(x)  \le (1 + 2 Y \beta(q)) \
\exp \left( - n^{q/(q+2)} \ [x/(K \delta)]^{2q/(q+2)} \right). \eqno(2.3)
$$

 {\bf Remark 2.1.} At the value $ q=1 $ we obtain the main result of article \cite{Lesign1};
at the value $ q \in (0,1) $  the main result of article \cite{Fan1}. \par

{\bf Example 2.2.} Assume now that

$$
T(\xi(i),x) \le C_1 \exp \left(-C_2(x/K)^q \ (\log(F(q,r) + x/K))^r \right), x > 0,
$$
where by definition
$$
F = F(q,r)=1, \ r \le 0; \ F(q,r)= \exp(q), \ r > 0.
$$

 Let us introduce the following vector - function $ L(q,r) = \{L(1;q,r),L(2;q,r)\}  $ of a variables
$(q,r):$
$$
L(1;q,r)= \frac{2q}{q+2}, \ L(2;q,r) = \frac{2r}{q+2}.
$$
 We conclude from theorem 2.1: $ \ T( n^{-1/2} \ \sum_i \xi(i), x) \le $

$$
\exp \left(-C_4 (x/K)^{L(1)} \ \left(\log(F(L(1),L(2))+x/K)^{L(2)} \right) \right) =: \exp \left( -G_{q,r}(x/K)  \right), \ x \le 2,
$$
where $ L(i) = L(i; q,r), \ i = 1,2; $   or equally

$$
Q_n(x) = T( n^{-1} \ \sum_i \xi(i), x) \le \exp \left( -G_{q,r}(x\sqrt{n}/K)  \right), \ x \ge 2.
$$

\vspace{3mm}

 We consider now moment, more exactly, $ L_p $ estimations for $ S(n)/n. $\par

\vspace{3mm}

{\bf Theorem 2.2.} Let $ \forall i \ \xi(i) \in L(p), \ p \ge 2; \ x \ge 1. $ Then

$$
Q_n(x)  \le x^{-p} \cdot (p-1)^p \cdot n^{-p/2} \cdot \left\{ \ n^{-1} \ \sum_{i=1}^n |\xi(i)|_p^2 \right\}^{p/2}. \eqno(2.4)
$$

 {\bf Proof.} The inequality

$$
\left|n^{-1/2} S(n) \right|_p \le (p-1) \left\{ \ n^{-1} \ \sum_{i=1}^n |\xi(i)|_p^2 \right\}^{1/2}. \eqno(2.5)
$$
is proved in  \cite{Ostrovsky0}; from (2.5) it follows

$$
\left|n^{-1} S(n) \right|_p \le (p-1)\cdot n^{-1/2} \cdot \left\{ \ n^{-1} \ \sum_{i=1}^n |\xi(i)|_p^2 \right\}^{1/2}. \eqno(2.6)
$$
 It remains to use the Markov's-Tchebychev's  inequality. \par

\vspace{3mm}

{\bf Remark 2.2.} The proposition of theorem 2.2 improves the correspondent moment estimate from the article
Lesign E., Volny D. \cite{Lesign1}. We write the exact value of constant in the right-hand side, which is less that
one in  \cite{Lesign1}; the non-improperness of the factor $ (p-1) $ is grounded also in  \cite{Ostrovsky0}.\par

\vspace{3mm}

{\bf Remark 2.3.} Assume that for some {\it interval} of a values $ p: \ p \in [2,a), $
where $ a = \const > 2  \ \forall i \ |\xi(i)|_p < \infty. $ It follows from the inequality (2.4)
then

$$
Q_n(x)  \le \inf_{ 2 \le p < a} \left[ x^{-p} \cdot (p-1)^p \cdot
n^{-p/2} \cdot \left\{ \ n^{-1} \ \sum_{i=1}^n |\xi(i)|_p^2 \right\}^{p/2} \right]. \eqno(2.7)
$$

\vspace{4mm}

\section{ Exactness of our estimates.} \par

\vspace{4mm}

 We prove in this section that the proposition of theorem 2.1 is non-improvable; we will construct an example
of martingales for which in (2.1) is attained asymptotical equality.\par
 Note that the exactness of proposition of theorem 2.2.is proved in \cite{Lesign1}.\par
We restrict ourselves only the case  when exists a positive value $ q $  such that
  for all the values $ n=1,2,\ldots $ and for all positive values $ x $

$$
\inf_i T(\xi(i),x)  \ge  \exp \left( - x^q \right), \
\sup_i T(\xi(i),x)  \le C_1 \ \exp \left( - x^q \right). \eqno(3.1)
$$

{\bf Theorem 3.1.} There exists a martingale $ (S(n), F(n))  $ with correspondent martingale-differences $ \xi(i) $
satisfying the conditions (3.1); and such that there exist finite positive values $ C_2, C_3 $  for which

$$
Q_n(1)  \ge C_2 \ \exp \left( - C_3 \ n^{q/(q+2)} \ \right).
$$

{\bf Proof.}   Assume at first $ q < 1; $ the case $ q=1 $ is investigated in \cite{Lesign1}.
Denote $ \beta = 1/(2+q), $ then $ 1-\beta = (q+1)/(2+q)  $ and $ \beta \in (0,1).  $ \par
We will follow the method offered by   Lesign E. and  Volny D. \cite{Lesign1}; namely, we intend to
built a very interest martingale of a view

$$
\xi(i) = \eta \cdot \zeta(i), \eqno(3.2)
$$
where  the r.v. $ \{\zeta(i) \}, \ \eta  $ are symmetrical distributed  and total independent:

$$
T(\eta,x)  = T(-\eta,x) =  \exp \left( - x^q \right), \eqno(3.3)
$$
the r.v. $ \zeta(i) $ are i.i.d., non-trivial and bounded: $ 0 < |\zeta(1)|_{\infty} < \infty. $ \par
 It is easy to verify that the r.v. $ \xi(i) = \eta \cdot \zeta(i) $ satisfy to the conditions (3.1).\par

  Further,

 $$
 Q_n(1) \ge {\bf P} \left(\sum_{i=1}^n \zeta(i) > C_4 \ n^{1-\beta} \right) \cdot
 {\bf P} \left(\eta > C_5 \ n^{\beta} \right) \stackrel{def}{=} Z_1 \cdot Z_2, \eqno(3.4)
 $$
where $ C_4, C_5 $ are suitable constants depending on $ q. $ We get:

$$
Z_2 \ge C_6 \ \exp \left( - C_7 \left[ n^{\beta} \right]^q \right) \ge
C_8 \ \exp \left( - C_9 \ n^{q/(q+2)} \right);  \eqno(3.5)
$$

$$
Z_1 = {\bf P} \left( \sum_{i=1}^n \zeta(i) > C_4 \ n^{(q+1)/(q+2)}  \right) =
{\bf P} \left( n^{-1/2}\sum_{i=1}^n \zeta(i) > C_4 \ n^{q/(2(q+2))}  \right).\eqno(3.6)
$$
 As long as here $ q \in (0,1), \ q/(2(q+2) \in (0,1/6) $ and we can apply the classical theory of great
deviations \cite{Petrov1}, \cite{Ibragimov1}, \cite{Linnik1}:

$$
Z_1 \ge C_{10} \ \exp \left( - C_{11} \ \left[ n^{q/(2(q+2))} \right]^2   \right) =
C_{10} \ \exp \left( - C_{11} \  n^{q/(q+2)}  \right). \eqno(3.7)
$$
 Multiplying the lower estimates for $  Z_1 $ and $  Z_2 $ we obtain what was required. \par

\vspace{3mm}

 Let now $ q > 1. $ Let us define the (unique) {\it integer} number $ s = s(q) $ such that
$ q-1 \le s < q; $ then

$$
\frac{s}{2(s+2)} \le  \frac{q}{2(q+2)}  < \frac{s+1}{2(s+3)}.  \eqno(3.8)
$$

 We will again construct a martingale differences by the formula (3.2), and the
r.v. $ \eta $  by the formula (3.3);  the r.v. $ \zeta(i) $ are again i.i.d., non-trivial, bounded
$ 0 < |\zeta(1)|_{\infty} < \infty. $ \par
 Moreover, we impose on the distribution $ \zeta(i) $ the following condition:

 $$
 {\bf E} \zeta^m(i)  = {\bf E} N^m(0,1), \ m = 2,3,4, \ldots, s+2, \eqno(3.9)
 $$
where $ N(0,1) $ denotes the standard normal distribution. \par
 We can use on the basis  of our assumptions  ((3.9) etc.) the theory of large deviations,
alike in the case $ 0 < q < 1; $ more exactly, we use theorem of Yu.V.Linnik  \cite{Linnik1}.\par
 After simple calculations we obtain what was required. \par

 \vspace{4mm}

\section{ Tail estimates for martingales under more strict conditions.} \par

\vspace{4mm}

 We impose in this section some strict condition on the martingale differences and obtain
an estimate for normed martingale at the same estimate as for sums of independent centered
variables. \par

 In order to formulate our result, we need to introduce some addition
notations and conditions. Let $ \phi = \phi(\lambda), \lambda \in (-\lambda_0, \lambda_0), \ \lambda_0 = \const \in (0, \infty] $
be some even strong convex which takes positive values for positive arguments twice continuous differentiable function 
in the intervals $ |\lambda| \ge 2, $ such that
$$
 \phi(0) = 0, \ \phi^{//}(0) \in(0,\infty), \ \lim_{\lambda \to \lambda_0} \phi(\lambda)/\lambda = \infty. \eqno(4.1)
$$
 We denote the set of all these function as $ \Phi; \ \Phi =
\{ \phi(\cdot) \}. $ \par
 We say that the {\it centered} random variable (r.v) $ \xi = \xi(\omega) $
belongs to the space $ B(\phi), $ if there exists some non-negative constant
$ \tau \ge 0 $ such that

$$
\forall \lambda \in (-\lambda_0, \lambda_0) \ \Rightarrow
{\bf E} \exp(\lambda \xi) \le \exp[ \phi(\lambda \ \tau) ].  \eqno(4.2)
$$
 The minimal value $ \tau $ satisfying (4.2) is called a $ B(\phi) \ $ norm
of the variable $ \xi, $ write
 $$
 ||\xi||B(\phi) = \inf \{ \tau, \ \tau > 0: \ \forall \lambda \ \Rightarrow
 {\bf E}\exp(\lambda \xi) \le \exp(\phi(\lambda \ \tau)) \}. \eqno(4.3)
 $$
 This spaces are very convenient for the investigation of the r.v. having a
exponential decreasing
 tail of distribution, for instance, for investigation of the limit theorem,
the exponential bounds of distribution for sums of random variables,
non-asymptotical properties, problem of continuous of random fields,
study of Central Limit Theorem in the Banach space etc.\par

  The space $ B(\phi) $ with respect to the norm $ || \cdot ||B(\phi) $ and ordinary operations is a Banach
space which is isomorphic to the subspace
consisted on all the centered variables of Orlicz's space $ (\Omega,F,{\bf P}), N(\cdot) $ with $ N \ - $ function

$$
N(u) = \exp(\phi^*(u)) - 1, \ \phi^*(u) = \sup_{\lambda} (\lambda u -
\phi(\lambda)).
$$
 The transform $ \phi \to \phi^* $ is called Young-Fenchel transform. The proof of considered assertion used
 the properties of saddle-point method and theorem of Fenchel-Moraux:
$$
\phi^{**} = \phi.
$$

 The next facts about the $ B(\phi) $ spaces are proved in \cite{Kozatchenko1}, \cite{Ostrovsky3}, p. 19-27:

$$
{\bf 1.} \ \xi \in B(\phi) \Leftrightarrow {\bf E } \xi = 0, \ {\bf and} \ \exists C = \const > 0,
$$

$$
T(\xi,x) \le \exp(-\phi^*(Cx)), x \ge 0.
$$
and this estimation is in general case asymptotically exact. \par

{\bf 2.}  Denote

 $$
 \overline{\phi}(\lambda) =  \sup_n \left[ n \phi(\lambda/\sqrt{n})  \right].  \eqno(4.4)
 $$

 Let $ \eta(i) $ be a sequence of i.i.d. random variables belonging to the space $ B(\phi). $ We reproduce
the following tail inequality for {\it normed} sum:

$$
T \left(n^{-1/2} \sum_{i=1}^n \eta(i),x \right) \le 2 \exp \left( - \{\overline{\phi}\}^*(x)  \right). \eqno(4.5)
$$

\vspace{3mm}

 Let  $ \phi(\cdot)  \in \Phi  $  and let $  \tilde{F} $ be any sub-sigma algebra of
 the source sigma-field $  F. $

\vspace{3mm}
{\bf Definition 4.1.}
\vspace{3mm}

 We say that the {\it centered} random variable (r.v) $ \xi = \xi(\omega) $
belongs to the space $ B(\tilde{F},\phi), $ if there exists some non-negative
{\it non-random}  constant $ \tau \ge 0 $ such that

$$
\forall \lambda \in (-\lambda_0, \lambda_0) \ \Rightarrow
{\bf E} \exp(\lambda \xi)/\tilde{F} \le \exp[ \phi(\lambda \ \tau) ].  \eqno(4.6)
$$
 The minimal value $ \tau $ satisfying (4.6) is said to be a $ B(\tilde{F},\phi) \ $ norm
of the variable $ \xi, $ write
 $$
 ||\xi||B(\tilde{F},\phi) = \inf \{ \tau, \ \tau > 0: \ \forall \lambda \ \Rightarrow
 {\bf E}\exp(\lambda \xi)/\tilde{F} \le \exp(\phi(\lambda \ \tau)) \}. \eqno(4.7)
 $$

  The space $ B(\tilde{F},\phi) $ with respect to the norm $ || \cdot ||B(\tilde{F},\phi) $ and ordinary
operations is also a Banach space. \par
 Our definition (4.1) is a direct generalization of a definition in  \cite{Liu1}, where is considered only the  case

$$
\phi(\lambda) = \phi_q(\lambda):= |\lambda|^q, \ |\lambda| > 1; \ \phi(\lambda) = \phi_q(\lambda)  = \lambda^2, \ |\lambda| \le 1,
$$
where $ \ q = \const \ge 1 $  (in our notations). \par

\vspace{3mm}

{\bf Theorem 4.1.} Suppose the martingale $ (S(n), F(n)) $ with correspondent martingale differences  $ \{\xi(i) \} $
is such that for some function $ \phi(\cdot) \in \Phi $

$$
\forall k = 1,2, \ldots \ \xi(k) \in B(F(k-1), \phi),
$$
and

$$
 \sup_k ||\xi(k)||B(F(k-1),\phi) \le 1 \ (\mod {\bf P}), \eqno(4.8)
$$
or equally

$$
\forall \lambda \in (-\lambda_0, \lambda_0) \ \Rightarrow
{\bf E}  \exp(\lambda \xi(k))/F(k-1) \le \exp[ \phi(\lambda) ] \ (\mod {\bf P}).  \eqno(4.9)
$$
 Then

$$
Q_n(x)  \le 2 \exp \left( - \{\overline{\phi}\}^*(x\sqrt{n})  \right). \eqno(4.10)
$$
 {\bf Proof.} Let $ \lambda \in (-\lambda_0, \lambda_0). $ We consider an exponential moment of the
martingale $ \lambda \ S(n). $ Indeed: $ {\bf E} \exp(\lambda S(n)) = $

$$
 {\bf E} \ {\bf E} \exp(\lambda S(n))/F(n-1) =
{\bf E} \left[ \lambda S(n-1) \ {\bf E} \exp(\lambda \xi(n))/F(n-1)   \right] \le
$$

 $$
{\bf E} (\lambda S(n-1)) \cdot \exp(\phi(\lambda)) \le {\bf E} ( \lambda S(n-2)) \cdot  \exp(2\phi(\lambda)) \ldots  \le
 $$

$$
\exp(n\phi(\lambda)),
$$
i.e. as in the case of summing of independent identical distributed random variables. \par
 It remains to use the formula (4.5), where instead $  x $ we write $ x \sqrt{n}. $ \par

\vspace{3mm}
{\bf Example 4.1.}  Suppose  in addition to the conditions of theorem 4.1 as in \cite{Liu1}
that
$$
\phi(\lambda) = \phi_q(\lambda), \ q = \const > 0,
$$
or equally alike in the article \cite{Liu1}

$$
\exists C \in (0,\infty), \ \vraisup_{\Omega} \  \sup_k {\bf E} \exp \left(C |\xi(k)|^q \right)/F(k-1) < \infty.
$$
 Define the variable  $ \gamma(q) $ as follows:

 $$
 \gamma(q) = \min(2,q), \ q \ge 2; \ \gamma(q) = 2q/(2+q), \ q \in (0,2).
 $$

We conclude on the basis of theorem 4.1:

$$
Q_n(x) \le C_1(q) \ \exp \left( -C_2(q) \ x^{\gamma(q)} \ n^{\gamma(q)/2 } \right).  \eqno(4.11)
$$
 Notice that the last estimate is non-improvable still for independent variables $  \xi(i). $ \par

{\bf Remark 4.1}  In the book \cite{Ostrovsky3}, p. 38-43 there are many of exponential
tail estimates  for sums of independent r.v. They may be generalize on the martingale case satisfying
the condition (4.8). \par
{\bf Remark 4.2} The non-asymptotical moment and tail estimates for heavy tailed polynomial martingales are
obtained in the article \cite{Ostrovsky7}. \par
 Notice that in the case of  martingales with heave tails the classical normalized sequence $ n^{-1/2} $ in
general case may be replaced.  For instance, if $ \xi(i) $ are i.i.d. symmetrical r.v.  with power tails of a view
(symmetrical Pareto distribution)
$$
T(\xi(i),x) = \min(1,x^{-r}), \ x > 0, \ r = \const \in (0,2),
$$
then the normalized sequence $  b(n) $ may be choose as follows: $ b(n) = n^{-1/r} $ (the Stable Limit Theorem).\par

 The case of {\it super-heavy} tails,  for example, distributions without variance, and when the r.v. $ \{ \xi(i) \} $
are independent centered summands was considered in \cite{Ostrovsky8} (non-asymptotical approach).

\vspace{4mm}

\section{ Inverse results to our estimates. }

\vspace{4mm}

 Suppose for any martingale  $ (S(n), F(n)) $ the following inequality is true:

$$
Q_n(x)  \le C_1 \  \exp \left( - C_2 \  n^{q/(q+2)} \ x^{2q/(q+2)} \right), \ x \ge 2, n \ge 1. \eqno(5.1)
$$
 Moreover, it is sufficient to assume the inequality (5.1) only for some {\it fixed} positive values $ x; $ say,
for the value $ x = 1: $

$$
Q_n(1) = {\bf P}( S(n) > n ) \le C_1 \  \exp \left( - C_2 \  n^{q/(q+2)} \right), \ n \ge 1. \eqno(5.2)
$$
 Question: give the possible  {\it lower} estimate for tail-function $ T(\xi(i), x).  $\par

{\bf Theorem 5.1.} Assume the martingale difference $ \{ \xi(i) \} $ are in addition independent
and identical distributed. Suppose also the condition (5.2) is satisfied. Then

$$
\sup_i T(\xi(i),x) \le C_3 \ \exp \left(-C_4 x^{2q/(q+2)} \right), \ x \ge 2.\eqno(5.3)
$$
{\bf Proof.} We will again follow the method offered by   Lesign E. and  Volny D. \cite{Lesign1}. Namely, we
will use the following fact:

$$
\underline{\lim}_{n \to \infty} \left[ \frac{Q_n(1)}{n {\bf P}(|\xi(1)| > 2n)} \right] \ge 1. \eqno(5.4)
$$
 It follows from (5.4) and (5.2) that for all sufficiently greatest values $ n, $  say $ n \ge 1 $

$$
\sup_i T(\xi(i),n) = T(\xi(1), n) \le C_5 \ n^{-1} \ \exp \left(-C_4 n^{2q/(q+2)} \right) \le
$$

$$
 C_5 \ \exp \left(-C_4 n^{2q/(q+2)} \right). \eqno(5.5)
$$

 This completes the proof of theorem 5.1. \par

\vspace{3mm}

{\bf Remark 5.1.} Assume that for some $ s = \const > 1 $

$$
Q_n(1) \le C \ n^{1-s}, \ n \ge 1.  \eqno(5.6)
$$
 We obtain by means of described method the tail inequality for $ \xi(i) $
in the case when the random variables $ \xi(i) $ are centered and i.i.d. that

 $$
\sup_i T(\xi(i),x)  \le C_5(s) \ x^{-s}, x \ge 1. \eqno(5.7)
 $$

 The last inequality implies that the r.v. $ C^{-1} \ \xi(i) $ belong to the unit
ball of the Lorentz space $ L_{s,\infty} = L_{s,\infty}(\Omega). $ Recall that the norm
in this space of the r.v. $  \xi $ may be defined as follows:

$$
||\xi||L_{s,\infty} = \sup_{x > 0} \left[x^s \ T(\xi,x) \right].
$$


\vspace{4mm}

\begin{thebibliography}{69}

\bibitem{Fan1}
 Fan X., Grama I, and Liu Q.
Large deviations for martingales with exponential condition.
arXiv:1111.1407 [math.PR] 6 Nov 2011.

\bibitem{Grama1}
 Grama I.G. On moderate deviations for martingales.
 Ann. Probab. Volume 25, Number 1 (1997), 152-183.

\bibitem{Grama2}
 Grama I. and Haeusler E., 2000. Large deviations for martingales via Cramer's
method. Stochastic Process. Appl. 85, 279–293.

\bibitem{Ibragimov1}
Ibragimov I. A. and Linnik Yu. V. (1965). Independent and Stationary Connected Variables.
Nauka, Moscow. (In Russian.) Mathematical Reviews (MathSciNet): MR202176
Zentralblatt MATH: 0154.42201

\bibitem{Kozatchenko1}
 Kozatchenko Yu. V., Ostrovsky E.I., The Banach Spaces of
random Variables of subgaussian type. Theory of Probab. and Math.
Stat. (in Russian). Kiev, KSU, {\bf 32,} 43 - 57 (1985).

\bibitem{Lesign1}
 Lesign E., Volny D. Large deviations for martingales.
Stochastic Processes and their Applications, {\bf 96,} 143 - 159 (2001)

\bibitem{Li1}
 Li Y, 2003. A martingale inequality and large deviations. Statist. Probab. Lett. 62,
317–321.

\bibitem{Linnik1}
Linnik Yu.V. Linit theorems for sums of independent r.v. taking into account of large
deviations. (1961), Theory Probab. Appl.,  {\bf 6},  $ N^o \ 2, $ p. 131-148;
 {\bf 6},  $ N^o \ 4, $ p. 363-360.

\bibitem{Liu1}
Liu Q., Watbled F. Exponential inequalities for martingales and asymptotic properties of the
free energy of directed polymers in random environment.
arXiv:0812.1719v1 [math.PR] 9 Dec 2008.

\bibitem{Ostrovsky0}
Ostrovsky E., Sirota L.  Moment and tail estimates for martingales and martingale transform,
with application to the martingale limit theorem in Banach spaces.
arXiv:1206.4964v1 [math.PR] 21 Jun 2012

\bibitem{Ostrovsky3}
 Ostrovsky E.I., 1999. Exponential estimations for Random Fields and
     its applications (in Russian). 1999, Obninsk, Russia, OINPE.

\bibitem{Ostrovsky4}
  Ostrovsky E.  Bide-side exponential and moment inequalities  for
      tails of distribution of Polynomial Martingales.
      arXiv: math.PR/0406532 v.1  Jun. 2004.

\bibitem{Ostrovsky7}
 Ostrovsky E. and  Sirota L. Moment and tail inequalities for polynomial martingales.
 The case of heavy tails.
arXiv:1112.2768v1 [math.PR] 13 Dez 2011

\bibitem{Ostrovsky8}
 Ostrovsky E. and  Sirota L. Non-improved uniform   tail estimates for normed sums of
independent random variables with heavy tails, with applications.
arXiv:1110.4879v1 [math.PR] 21 Oct 2011.

\bibitem{Saulis1}
Saulis, L. and Statulevicius, V. (1989). Limit Theorems for Large Deviations. Mokslas, Vilnius. (In Russian.)
Mathematical Reviews (MathSciNet): MR93e:60055b

\bibitem{Petrov1}
 Petrov V.V. Limit theorems of probability theory, volume 4 of Oxford Studies in Probability.
The Clarendon Press Oxford University Press, New York, 1995. Sequences of independent random
variables, Oxford Science Publications.

\bibitem{Volny1}
Volny D. Approximating martingales and the central limit theorem for strictly  stationary
processes. Stoch. processes ant their applic., {\bf 44}, (1993). 41-74.

\bibitem{Ratchkauskas1}
Ratchkauskas A. Large deviations for martingales with some applications.
Acta Applicandae Mathematicae, Volume 38, Number 1 (1995), 109-129, DOI: 10.1007/BF00992617


\end{thebibliography}
\end{document}